\documentclass[11pt]{article}
\usepackage{amssymb,latexsym, amsmath}
\usepackage{amscd, amsthm}

\newcommand{\be}{\begin{eqnarray}}
\newcommand{\ee}{\end{eqnarray}}
\newcommand{\bew}{\begin{eqnarray*}}
\newcommand{\eew}{\end{eqnarray*}}
\makeatletter \@addtoreset{figure}{section}
\def\thefigure{\thesection.\@arabic\c@figure}
\def\fps@figure{h, t}
\@addtoreset{table}{bsection}
\def\thetable{\thesection.\@arabic\c@table}
\def\fps@table{h, t}
\@addtoreset{equation}{section}

\newif\ifamsfonts
\amsfontstrue \ifamsfonts
\font\twlbbb=msbm10 scaled\magstep1 \font\egtbbb=msbm8
\font\sixbbb=msbm6
\newfam\mathbbfam
\textfont\mathbbfam=\twlbbb \scriptfont\mathbbfam=\egtbbb
\scriptscriptfont\mathbbfam=\sixbbb

\makeatother

\newtheorem{corollary}{Corollary}[section]
\newtheorem{theorem}{Theorem}[section]
\newtheorem{proposition}[theorem]{Proposition}

\newfont{\tenbi}{cmbxti10}

\begin{document}
\title{Control of mechanical systems on Lie groups and ideal hydrodynamics}

\author{\Large Mikhail V. Deryabin \footnote{Mads Clausen
Instituttet, Syddansk Universitet, Grundtvigs All{\'e} 150, DK-6400
S{\o}nderborg. Email:~mikhail@mci.sdu.dk }}

\date{}

\maketitle

\small

\begin{center}
\textbf{Abstract}
\end{center}

In contrast to the Euler-Poincar{\'e} reduction of geodesic flows of
left- or right-invariant metrics on Lie groups to the corresponding
Lie algebra (or its dual), one can consider the reduction of the
geodesic flows to the group itself. The reduced vector field has a
remarkable hydrodynamic interpretation: it is a velocity field for a
stationary flow of an ideal fluid. Right- or left-invariant symmetry
fields of the reduced field define vortex manifolds for such flows.

Consider now a mechanical system, whose configuration space is a Lie
group and whose Lagrangian is invariant to left translations on that
group, and assume that the mass geometry of the system may change
under the action of internal control forces. Such system can also be
reduced to the Lie group. With no controls, this mechanical system
describes a geodesic flow of the left-invariant metric, given by the
Lagrangian, and thus its reduced flow is a stationary ideal fluid
flow on the Lie group. The standard control problem for such system
is to find the conditions, under which the system can be brought
from any initial position in the configuration space to another
preassigned position by changing its mass geometry. We show that
under these conditions, by changing the mass geometry, one can also
bring one vortex manifold to any other preassigned vortex manifold.

\large


\section {Introduction}


For the Euler top, the Hamiltonian vector field on the cotangent
bundle $T^* SO(3)$ can be uniquely projected onto the Lie algebra
$\mathfrak{so}(3)$ -- this is a classical reduction, known in the
general case as the Euler-Poincar{\'e} reduction. In the 30's,
E.T.Whittaker suggested an ''alternative'' reduction procedure for
the Euler top: by fixing values of the Noether integrals, the
Hamiltonian vector field can be uniquely projected from $T^* SO(3)$
onto the group $SO(3)$ \cite{Whi}. The Whittaker reduction is valid
for any Hamiltonian system on a cotangent bundle $T^* G$ to a Lie
group $G$, provided the Hamiltonian is invariant under the left (or
right) shifts on the group $G$. An important example of such
Hamiltonian systems is a geodesic flow of a left-(right-)invariant
metric on a Lie group.

If we reduce a Hamiltonian system to the Lie group $G$, and then
factorize the reduced vector field by the orbits of its symmetry
fields, then, by the Marsden-Weinstein theorem, we get the same
Hamiltonian system on a coadjoint orbit on the dual algebra
$\frak{g}^*$, as if we first reduced the system to the dual algebra
$\frak{g}^*$, and then to the coadjoint orbit (see also \cite{Arn},
Appendix 5). Thus, the Whittaker reduction can be regarded as a part
of the Marsden-Weinstein reduction of Hamiltonian systems with
symmetries \cite{MW}.

In contrast to the Marsden-Weinstein reduction, it has not been
payed much attention to the Whittaker reduction alone. However, it
is itself worth studying. It turns out that a vector field, reduced
to a Lie group $G$ has a remarkable {\it hydrodynamic
interpretation}: it is a velocity field for a stationary flow of an
ideal fluid, that flows on the group $G$ (viewed as a Riemannian
manifold), and is incompressible with respect to some left-(or
right-)invariant measure on $G$, see \cite{Der06, Koz1, Koz2, Koz_b}
for details. The reduction to a Lie group is also useful for a
series of applications, which include stability theory,
noncommutative integration of Hamiltonian systems, discretization,
differential geometry of diffeomorphism groups and see, e.g.,
\cite{Koz_b, Koz3, DF, Der06, Fed_1}.

In this article, we first review the Whittaker reduction and its
hydrodynamic essence, and provide an explicit expression for the
reduction of a geodesic flow of a left- or right-invariant metric
onto a Lie group. For any Lie group we find both the reduced vector
field and its ''symmetry fields'', i.e., left- or right-invariant
fields on the group that commute with our reduced vector field.
These fields have also a hydrodynamical meaning: these are the {\it
vortex} vector fields for our stationary flow (i.e., they annihilate
the vorticity 2-form), cf. \cite{Koz_b}, \cite{Der06}. The
distribution of the vortex vector fields in always integrable, thus
they define a manifold, that we call the {\it vortex manifold}.
Typically, these manifolds are tori.

Next, we consider the following control problem. We study mechanical
systems, whose configuration space is a Lie group and whose
Lagrangian is invariant to left translations on that group, and we
assume that the mass geometry of the system may change under the
action of internal control forces. Such systems can also be reduced
to the Lie group, and they also have an interesting hydrodynamic
interpretation: the reduced vector field is the velocity of a
stationary flow of an electron gas (with no controls, this
mechanical system describes a geodesic flow of the left-invariant
metric, given by the Lagrangian, and thus its reduced flow is a
stationary ideal fluid flow). Notice that without relating to
hydrodynamics, controlled systems on Lie groups were studied in many
works, see, e.g., \cite{A} and references therein.

The standard control problem for such systems is to find the
conditions, under which the system can be brought from any initial
position on the Lie group to another preassigned position by
changing its mass geometry. We show that under these conditions, by
changing the mass geometry, one can bring the whole vortex manifold
to any other preassigned vortex manifold.

As an example, we consider the $n$-dimensional Euler top. We write
down the reduced controlled system explicitly, find the vortex
manifolds, which typically (when the momentum matrix has the maximal
rank) are tori, and show that, by changing the mass geometry, every
such vortex manifold can be transformed to any other vortex
manifold.

In the Appendix we study the Whittaker reduction for nonholonomic
systems, and formulate and discuss the standard controllability
conditions.

\section{Reduction of a geodesic flow to a Lie group} \label{Sec_2}

We start with some basic facts on coadjoint representations, inertia
operators on Lie algebras and the Euler equations (see, e.g., \cite
{AK}). Let $G$ be an arbitrary Lie group, $ \mathfrak {g} $ be its
Lie algebra, and $\mathfrak {g}^*$ be the corresponding dual
algebra. The group $G$ may be infinite-dimensional, and not
necessarily a Banach manifold, but we assume that the exponential
map $\exp: \frak{g} \to G$ exists.

Any vector $\dot {g} \in T_g G$ and any covector $m\in T^*_g G$ can
be translated to the group unity by the left or the right shifts. As
the result we obtain the vectors $ {\omega}_c, {\omega}_s \in
\mathfrak{g}$ and the momenta $m_c, m_s\in \mathfrak {g}^ *$:
$$
\omega_c = L_{g^{-1}*} \dot{g}, \quad \omega_s = R_{g^{-1}*}
\dot{g}, \quad m_c = L^*_g m, \quad m_s = R^*_g m.
$$
The following relation plays the central role in the sequel:
\begin{equation} \label{mcs}
m_c= Ad^*_ g m_s,
\end{equation}
$Ad ^*_g: \mathfrak {g}^*\to\mathfrak {g}^*$ being the group
coadjoint operator. Let us fix the ''momentum in space'' $m_s$. Then
relation (\ref{mcs}) defines a coadjoint orbit. The Casimir
functions are the functions of the ''momentum in the body'' $m_c$,
that are invariants of coadjoint orbits. For example, for the Euler
top, the Casimir function is the length of the kinetic momentum.

Let $A: \mathfrak {g} \to \mathfrak {g}^*$ be a positive definite
symmetric operator ({\it inertia operator}) defining a scalar
product on the Lie algebra. This operator defines a left- or
right-invariant inertia operator $A_g$ (and thus a left- or
right-invariant metric) on the group $G$. For example, in the
left-invariant case, $A_g = {L_{g}^*}^{-1} \, A \, L_{g^{-1} *}$.

Let the metric be left-invariant. The geodesics of this metric are
described by the Euler equations
 \be \label{geod}
\dot {m_c} = ad ^ * _ {A ^ {-1} m_c} m_c,
 \ee
Here $ad^*_ {\xi}: \mathfrak {g}^* \to \mathfrak {g}^*$ is the
coadjoint representation of $\xi \in \mathfrak{g}$. Given a solution
of the Euler equations $\omega_c = A^{-1} m_c$, the trajectory on
the group is determined by the relation
 \be \label{Pois:eq}
L_{g^{-1}*}\dot{g}={\omega}_c.
 \ee
The Euler equations follow from the fact that "the momentum in
space" $m_s $ is constant, whereas "the momentum in the body" $m_c$
is obtained from $m_s $ by (\ref{mcs}), see \cite{AK}.

\noindent{\bf Remark.} Strictly speaking, in the infinite
dimensional case the operator $A$ is invertible only on a regular
part of the dual algebra $\mathfrak {g}^*$. In our case this means,
that some natural restriction on values of $m_s $ (or $m_c $) have
to be imposed (see~\cite {AK}).

The Euler equations can be considered as Hamilton's equations on the
dual algebra, where the Hamiltonian equals $H = \frac{1}{2}
(A^{-1}m, m)$, $m \in \mathfrak{g}^*$, and the Poisson structure is
defined by the following Poisson brackets. For two functions $F(m)$
and $G(m)$ on the dual algebra $\mathfrak{g}^*$,
 \bew
\{F,G\} = \left(m, \left[dF(m), dG(m)\right] \right),
 \eew
where $dF(m), dG(m) \in \mathfrak{g}$ are the differentials of
functions $F$ and $G$, and $[\xi, \eta] = ad_\xi \eta$ is the
commutator (adjoint action) on the Lie algebra~$\mathfrak{g}$.

Let now
 \bew
H = \frac{1}{2} (A^{-1}m_c, m_c)  + (\lambda, m_c),
 \eew
where $\lambda \in \mathfrak{g}$ is a constant vector. Then Equation
(\ref{geod}) becomes
 \be \label{Eu:eqs}
\dot {m_c} = ad ^ * _ {A^{-1} m_c + \lambda} m_c,
 \ee
and the velocity $\omega_c = A^{-1} m_c + \lambda$

In case of a right-invariant metric $m_c$ is constant, the Euler
equations read $\dot {m_s}=-ad^*_{A^{-1} m_s + \lambda}m_s$, and the
trajectory on the group is determined by the equation $R_{g^{-1} *}
\dot {g} = {\omega}_s$.

The result of the reduction onto the group is a vector field $v (g)
\in TG$ such that the trajectory on the group is defined by the
equation $\dot{g}=v (g)$. The field $v(g)$ will be referred to as
{\it reduced}.

\begin{proposition} \label{Pr1} {\bf (The Whittaker reduction)}
For the case of the left-invariant or the right-invariant metric,
the vector field $v (g) $ has the form
\begin {equation} \label{v_lev}
v (g) =L _ {g *} (A ^ {-1} Ad ^ * _ g m_s + \lambda)
\end{equation}
and, respectively
\begin{equation} \label{v_prav}
v (g) =R _ {g *} (A ^ {-1} Ad ^ * _ {g ^ {-1}} m_c + \lambda).
\end{equation}
Here $m_s$, respectively $m_c$, is constant.
\end{proposition}

Notice that in Proposition \ref{Pr1}, to find the reduced vector
field we do not need the Hamiltonian equations on $T^*G$ and the
explicit expression for the Noether integrals. We only need the Lie
group structure and the inertial operator. This is important for
generalizations to the infinite-dimensional case. Unlike for the
Marsden-Weinstein reduction, we do not have to assume nondegeneracy
conditions on the momenta $m_s$ or $m_c$.
\medskip

\noindent{ \it Proof.} We consider only the case of the
left-invariant metric; for the right-invariant case the proof is
similar. Relation (\ref{mcs}) determines the function $m_c=m_c (m_s,
g)$ on the group $G$ depending on $m_s$ as a parameter. From the
equality $\omega_c =A ^ {-1} m_c + \lambda$ and $L _ {g ^ {-1} *}
\dot {g} = {\omega} _c $ follows that for any $g \in G $,
$$
L _ {g ^ {-1} *} \dot {g} = A ^ {-1} m_c (m_s, g) + \lambda
$$
which implies (\ref{v_lev}). $\Box$
\medskip

In Appendix A we consider the case, when the inertia operator is not
left- or right-invariant, i.e, $A = A(g)$. Some nonholonomic systems
have this form. It turns out that system of equations
(\ref{Pois:eq}), (\ref{Eu:eqs}) can still be reduced to the group
$G$ (although now Equation (\ref{Eu:eqs}) cannot be separated).

Even if $\lambda = 0$, the reduced vector fields (\ref{v_lev}),
(\ref{v_prav}) are is in general neither left- nor right-invariant.
An important exception is when the inertia operator defines a
Killing metric on the Lie algebra. However, the reduced {\it
covector} fields are always right- or left-invariant.

\begin{proposition} \label{Pr_LR}
Let $\lambda = 0$. If the metric is left-invariant, then the reduced
covector field $m(g) = A_g v(g)$ is right-invariant.
\end{proposition}

{\it Proof.}
$$
m(g) = {L_{g}^*}^{-1} A L_{g^{-1} *} v(g) = {L_{g}^*}^{-1} A A^{-1}
L^*_g {R^*_{g}}^{-1} m_s = {R^*_{g}}^{-1} m_s.
$$
$\Box$

Let $w (g) \in TG $ be a  right-invariant vector field on the group
$G$, which is defined by a vector $\xi \in \mathfrak {g}$: $w (g) =
R _ {g *} \xi $.

We fix a momentum $m_s$.

\begin{theorem} \label{Th1}
For the momentum $m_s$ fixed, the vector field $w (g) $ on $G$ is a
symmetry field of the reduced system $v (g)$ if and only if the
vector $\xi$ satisfies the condition
\begin{equation} \label{main_cond}
ad _ {\xi} ^ * m_s = 0.
\end{equation}
\end{theorem}

In the finite-dimensional case this means that the flows of the
vector fields $v(g), w(g)$ on the group commute. In the
infinite-dimensional case one should be more accurate: the equation
$\dot {g} = v (g)$ is a partial integral-differential equation,
rather than an ordinary differential equation, hence, strictly
speaking, it is not clear if it has a solution. On the other hand,
equation $\dot {g} = R_{g *} \xi$ always has a solution, which is a
one-parametric family of the left shifts on the group $G$: $g \to (
\exp {\tau \xi}) g$, see, for example, ~\cite{Che}, as we have
assumed that the exponential map exists. Notice also that, in view
of Proposition \ref{Pr_LR}, under the assumption of Theorem
\ref{Th1}, the Lie derivative $L_{w(g)} m(g) = 0$.
\medskip

\noindent{\it Proof of Theorem \ref{Th1}}. The vector fields $w(g)$
and $L_{g*} \lambda$ commute, as right-invariant fields always
commute with left-invariant fields. Thus, it is sufficient to show
that
$$
v ((\exp {\tau \xi}) g) = L _ {(\exp {\tau \xi}) *} v (g)
$$
if and only if the condition of the theorem is fulfilled. Indeed,
\begin{align*}
v ((\exp {\tau \xi}) g) & =L _ {(\exp {\tau \xi}) g *} A ^ {-1}
Ad ^ * _ {(\exp {\tau \xi}) g} m_s \\
& = L_ {\exp {\tau \xi} *} L _ {g *} A ^ {-1} Ad ^ * _ g (Ad ^ * _
{\exp {\tau \xi}} m_s).
\end{align*}
The last term equals $L_ {(\exp {\tau \xi}) *} v (g)$ for any $g\in
G$ if and only if
$$
Ad^*_ {\exp {\tau \xi}} m_s = m_s
$$
for all values of the parameter $\tau$. Differentiating the last
relation by $\tau$ we arrive at the statement of the theorem. $\Box$

\section{Stationary flows on Lie groups} \label{Sec:hydro}

We now formulate some results on the hydrodynamics character of the
reduced vector fields from the previous section. Consider first the
Euler equations for an ideal incompressible fluid, that flows on a
Riemannian manifold $M$:
$$
\frac{\partial v}{\partial t} + \nabla_v v = -\nabla p, \quad
\mbox{div} \, v = 0,
$$
where $\nabla_v v$ is the covariant derivative of the fluid velocity
vector $v$ by itself with respect to the Riemannian connection and
$p$ is a pressure function.

Consider a geodesic vector field $u$ on the manifold $M$. Locally it
always exists, but it may not be defined globally on $M$ -- take a
two-sphere as a simple example. Then $u$ is a stationary flow of the
ideal fluid with a constant pressure. Indeed, as $u$ is a geodesic
vector field, its derivative along itself is zero: $\nabla_u u = 0$.

{\bf Remark.} The converse is of course not true: there are
stationary flows that are not geodesics of the Riemannian metric.

The stationary flows with constant pressure form a background for
hydrodynamics of Euler equations on Lie groups. Consider a
Hamiltonian system on a finite-dimensional Lie group $G$, with a
left-invariant Hamiltonian, which is quadratic in the momenta (in
terms of Section~\ref{Sec_2}, vector $\lambda = 0$). This
Hamiltonian defines a left-invariant metric on the Lie group $G$. As
we reduce this system to the group, the reduced vector field is
globally defined on $G$, and is a geodesic vector field of the
Riemannian metric, defined by the left-invariant Hamiltonian, and it
defines a stationary flow of an ideal fluid on $G$.

Thus, the reduced vector field (\ref{v_lev}) (and (\ref{v_prav})) is
the velocity vector field for a stationary flow on the Lie group $G$
with left- (right-) invariant metric. An immediate corollary of
Proposition \ref{Pr_LR} is

\begin{proposition} \label{Pr3.1}
There is an isomorphism between the stationary flows with constant
pressure, defined by a left-invariant metric on a finite-dimensional
Lie group $G$, and the space of right-invariant covector fields on
this group.
\end{proposition}

{\bf Remark.} Stationary flows with constant pressure play an
important role in studying the differential geometry of
diffeomorphism groups, see \cite{BR, KM, M}: they define asymptotic
directions on the subgroup of the volume-preserving diffeomorphisms
of the group of all diffeomorphism. Proposition \ref{Pr3.1} is a
generalization of \cite{Pa}, where it was shown that every
left-invariant vector field on a compact Lie group equipped with a
bi-invariant metric is asymptotic: if a Hamiltonian defines the
bi-invariant metric on the Lie algebra, then the reduced vector
field (\ref{v_lev}) is itself left-invariant. Moreover, its flow
(which are right shifts on the Lie group $G$) are isometries of this
metric (see, e.g., \cite{DNF}).

Recall now that the reduced covector field is right-invariant
(Proposition \ref{Pr_LR}). Thus, the condition $ad^*_\eta m_s = 0$
is equivalent to $L_{\eta(g)} m(g) = 0$, where $m(g)$ is the
right-invariant 1-form (being equal to $m_s$ at $g = id$), and
$\eta(g) = R_{*g} \eta$ is the right-invariant symmetry field. By
the homotopy formula,
$$
0 = L_{\eta(g)} m(g) = i_{\eta(g)} dm(g) + d ({\eta(g)}, m(g)) =
i_{\eta(g)} dm(g),
$$
as $(\eta(g), m(g)) = (\eta, m_s) = const$ for all $g$ (both vector
and covector fields are right-invariant).

We now define a {\it vortex vector field}, as an annihilator of the
vorticity 2-form. Then the condition $i_{\eta(g)} dm(g) = 0$ is
exactly the definition of a vortex field. Thus, we have proved

\begin{proposition}
Any symmetry field to the reduced vector field is a vortex vector
field.
\end{proposition}

Vortex vectors, i.e., vectors $\xi \in \mathfrak{g}$ that satisfy
condition (\ref{main_cond}), are the isotropy vectors. We now review
some classical results on the isotropy vectors and the Casimir
functions, see, e.g., \cite{AG} for details, and adapt them to our
case.

\begin{proposition} \label{Pr01}
The distribution of the isotropy vectors in integrable.
\end{proposition}

The Proposition says that if vectors $\xi_1, \xi_2 \in \mathfrak{g}$
satisfy condition (\ref{main_cond}), then the vector $[\xi_1, \xi_2]
= ad_{\xi_1} \xi_2$ also satisfies this condition, which is a simple
consequence of the Jacobi identity. The integrable distribution of
the isotropy vectors defines a manifold (at least locally), that we,
following \cite{Koz_b}, call a {\it vortex manifold}.

The isotropy vectors $\xi \in \mathfrak{g}$ form a Lie subalgebra
$\mathfrak{h} \subset \mathfrak{g}$, called an isotropy algebra for
the coadjoint orbit $m = Ad^*_g m_s$. If the differentials of the
Casimir functions form a basis of the isotropy algebra
$\mathfrak{h}$, then $\mathfrak{h}$ is Abelian. In general, an
isotropy algebra is not necessarily Abelian. A very simple example
is $G = SO(3)$: if the ''momentum in space'' $m_s = 0$, then
$\mathfrak{h} = \mathfrak{so}(3)$. However, in the
finite-dimensional case isotropy algebras are Abelian on an open and
dense set in $\mathfrak{g}^*$ (the Duflo theorem). Thus, the
corresponding vortex manifolds (that pass through the group unity)
are commutative subgroups of the Lie group~$G$. Notice that in the
infinite-dimensional case, vortex fields cam still define a certain
commutative subgroup, which can also be referred to as ''vortex
manifold''.

Vortex manifolds have always the dimension of the same parity as the
Lie group dimension. This is a simple corollary of the fact that
coadjoint orbits are always even-dimensional (also the degenerate
ones), see, e.g., \cite{AG}.

One can show that if the Hamiltonian has also terms, linear in the
momenta (in the other words, if $\lambda \ne 0$), then the reduced
field has the following hydrodynamic sense: it is the velocity of
the stationary flow for the electron gas, which satisfies an
''infinite conductivity equation'', again, with a constant pressure,
see~\cite{AK}.

\section{Control on Lie groups and vortex manifolds}

Consider a Lagrangian system on a tangent bundle $TG$ to a Lie
algebra $G$, with the Lagrangian, which is left-invariant under the
action of the Lie group $G$.

In order to introduce the controls in our system, we consider
Lagrangians on $TG$ of the following form:
$$
L(\omega,u) = \frac{1}{2} \left( A(\omega + \sum_{i = 1}^k u_i
\lambda_i), \omega + \sum_{i = 1}^k u_i \lambda_i \right),
$$
where $\omega \in \mathfrak{g}$ is the system velocity, $\lambda_i
\in \mathfrak{g}$ are constant vectors, $u_i(t) \in \mathbb{R}$ are
controls, and $A: \mathfrak{g} \to \mathfrak{g}^*$ is the inertia
operator. Notice that the dimension $k$ of the control vector $u(t)$
may be lower than the dimension of the Lie algebra. We assume that
there is a positive constant $\epsilon$, such that $\|u(t)\| \le
\epsilon$, i.e., our controls are always bounded. Physically, these
controls mean that we can change the system mass geometry by
internal forces.

The Euler equations (\ref{Eu:eqs}) are:
$$
\dot{m} = ad^*_{\omega} m,
$$
where the momentum $m = A(\omega + \sum_{i = 1}^k u_i \lambda_i) \in
\mathfrak{g}^*$. The system, reduced to the group $G$, is (cf.
(\ref{v_lev}):
\begin{equation} \label{red_control}
\dot{g} = L_{g *} \left( A^{-1} Ad^*_g m_s - \sum_{i = 1}^k u_i
\lambda_i \right) = v_\lambda (g).
\end{equation}

From Theorem \ref{Th1} and Proposition \ref{Pr01} follows the
following result.

Suppose that System (\ref{red_control}) is controllable (we
formulate corresponding conditions in the Appendix \ref{Sec:AppB}),
and we assume that the controls $u(t)$ are piecewise constant
functions. We fix the ''momentum in space'' $m_s$: with $m_s$ fixed,
so are the vortex manifolds.

\begin{theorem} \label{Th_vort}
By applying controls $u(t)$, one can transform any vortex manifold
$H_1$ to any other prescribed vortex manifold $H_2$, such that the
following diagram is commutative:
 \be
\begin{array}{ccc}
H_1  & \stackrel{g_{v_\lambda}^t}{\to} & H_2\\
g^s_w \downarrow &  & \downarrow g^s_w \\
H_1 & \stackrel{g_{v_\lambda}^t}{\rightarrow} & H_2,
\end{array}
 \ee
where by $w$ we denote vortex vector fields for the given momentum
$m_s$, $g^s_w$ being its phase flow, and $g^t_{v_\lambda}$ is the
phase flow of System (\ref{red_control}).
\end{theorem}

This theorem is a reflection of a well-known fact that vortex lines
are frozen into the flow of an ideal fluid.

{\it Proof.} By Theorem \ref{Th1}, the vector fields $v_\lambda(g)$
and the vortex fields $w(g)$ commute (the vortex fields are
right-invariant, while the vectors $L_{g*} \lambda_i$ are
left-invariant, and we have also assumed that $u(t)$ is piecewise
constant). Pick up the controls (i.e., functions $u(t)$), that send
a point $h_1 \in H_1$ to a point $h_2 \in H_2$. Then the same
controls send a point $g^s_w h_1$ to $g^s_w h_2$, due to
commutativity, which proves the theorem. $\Box$

A simple corollary is that all vortex manifolds, that correspond to
the same value of the momentum $m_s$, are homotopic to each other.
Another observation is that an electron gas, flowing on a Lie group,
can be controlled by changing an external electro-magnetic field.

As an example, we consider the control problem for an
$n$-dimensional rigid body with a fixed point in ${\mathbb R}^{n}$
($n$-dimensional top). We follow the reduction procedure, suggested
in \cite{DF}.

Let $\mathfrak {so}(n)$ be the Lie algebra of $SO(n)$, $R \in SO(n)$
be the  rotation matrix of the top, $\Omega_c=R^{-1} \dot R
\in\mathfrak {so}(n)$ be its angular velocity in the moving axes,
and $M_c\in \mathfrak {so}^* (n)$ be its angular momentum with
respect to the fixed point of the top, which is also represented in
the moving axes.

The angular momentum in space $M_s=Ad^*_{R^{-1}} M_c \equiv R M_c
R^{-1}$ is a constant matrix, and the Euler equations have the
following matrix form generalizing the classical Euler equations of
the rigid body dynamics
\begin{equation} \label {2.8}
\dot {M} _c + [\Omega_c, M_c] =0.
\end{equation}

We assume, that the inertia operator of $A\, : \mathfrak {so} (n)
\to \mathfrak {so}^*(n)$ is defined by the relation $\Omega_c=A ^
{-1} M = UM+MU $, where $U$ is any constant nondegenerate operator.
Thus the system (\ref {2.8})  is a closed system of $n (n-1) /2$
equations, which was first written in an explicit form  by F. Frahm
(1874) \cite{FRahm}. As was shown in \cite {Man} (for $n=4$, in
\cite {FRahm}), with the above choice of the inertia tensor, the
system ( \ref {2.8}) is a completely integrable Hamiltonian system
on the coadjoint orbits of the group $SO (n)$ in $\mathfrak
{so}^*(n)$.

Now we fix the angular momentum $M_s$ (and, therefore, the coadjoint
orbit) and assume that rank $M_s = k\le n$ ($k$ is even). Then,
according to the Darboux theorem (see, e.g., \cite{AG}), there exist
$k$ mutually orthogonal and fixed in space vectors $x^{(l)},
y^{(l)}$, $l=1, \dots, k/2$ such, that $|x ^ {(l)}|^2 =
|y^{(l)}|^2=h_l$, $h_l =$const, and the momentum can be represented
in the form
\begin{equation} \label {r5.87}
M_s =\sum _ {l=1} ^ {k/2} x ^ {(l)} \wedge y ^ {(l)}, \quad \mbox
{that  is} \quad M_s = {\cal X} ^T {\cal Y} - {\cal Y} ^T {\cal X},
\end{equation}
where ${\cal X}^T= (x^{(1)} \, \cdots \, x^{(k/2)}), \quad {\cal
Y}^T = (y^{(1)} \, \cdots \, y ^ {(k/2)})$, $x^{(l)}\wedge y^{(l)}
=x ^ {(l)} \otimes y^{(l)} - y^{(l)}\otimes x^{(l)}$, and $(\,)^T$
denotes transposition. Under these conditions on $x^{(l)}, y^{(l)}$
the set of $k\times n$ matrices ${\cal Z}= (x ^{(1)} \, y^{(1)} \,
\cdots \, x^{(k/2)} \, y^{(k/2)})^T$ forms the Stiefel variety
${\cal V} (k, n)$ (see, for example, \cite {DNF}).

The momentum in the body $M_c$ has the same expression as (\ref
{r5.87}), but here the components of matrices $\cal{X}, \cal{Y}$ are
taken in a frame attached to the body, see (\ref{mcs}).

Since the above vectors are fixed in space, in the moving frame they
satisfy the Poisson equations, which are equivalent to matrix
equations
\begin{equation} \label {r5.88}
\dot {\cal X} = {\cal X} \Omega_c, \quad \dot {\cal Y} = {\cal Y}
\Omega_c.
\end{equation}

Now we set $\Omega_c=UM_c+M_cU $ and substitute this expression into
(\ref{r5.88}). Then taking into account (\ref {r5.87}), we obtain
the following dynamical system on ${\cal V} (k, n)$
\begin{eqnarray} \label {r588}
\dot {\cal X}&={\cal X} [U ({\cal X} ^T {\cal Y} - {\cal Y} ^T {\cal
X})
+ {\cal X} ^T {\cal Y} U], \nonumber \\
\dot {\cal Y}&={\cal Y} [U ({\cal X} ^T {\cal Y} - {\cal Y} ^T {\cal
X}) - {\cal Y} ^T {\cal X} U].
\end{eqnarray}

Notice that in the case of maximal rank $k $ ($k=n$ or $k=n-1$), the
Stiefel variety is isomorphic to the group $SO(n)$, and the
components of vectors $x^{(1)} \, y^{(1)} \, \cdots \,x^{(k/2)} \,
y^{(k/2)} $ form redundant coordinates on it. Thus the system
(\ref{r588}) describes required reduced flow (\ref {v_lev}) on
$SO(n)$.

The representation (\ref {r5.87}) is not unique: rotations in
2-planes spanned by the vectors $x^{(l)}, y^{(l)}$ in ${\mathbb R}
^n$ (and only they), leave the angular momentum $M$ invariant (in
the case of the maximal rank). As a result, the system (\ref {r588})
on $SO(n)$ has $k/2$ vortex vector fields $w_1(g), \dots, w_{k/2}
(g)$, which are generated by the right shifts of vectors $\xi^ {l}
\in\mathfrak {so}(n)$, such that $ad_{\xi^{l}}^ * M_s \equiv
[\xi^{l},M]=0$, cf. Section \ref{Sec:hydro}. In the redundant
coordinates the fields take the form
$$
\dot x ^ {(l)} = (x ^ {(l)}, x ^ {(l)}) y ^ {(l)}, \quad \dot {y ^
{(l)}} = - ( y ^ {(l)}, y ^ {(l)}) x ^ {(l)}, \quad l=1, \dots, k/2.
$$

One can easily see that in the case of maximal rank of the momentum
matrix, the corresponding vortex manifolds are $k/2$-{\it
dimensional tori}. This is a general fact: if a Lie group is
compact, then the vortex manifolds are compact manifolds, and, by
the Duflo theorem, for a dense set of the momenta $m_s$, the vortex
manifolds are tori (in our case, this dense set is determined by the
condition that the momentum rank is maximal). The torus, that passes
through the group unity, is called the {\it maximal torus for the
Lie group}; maximal tori play an important role in classification of
compact Lie groups.

One can furthermore show that if the rank of the momentum is not
maximal, the vortex manifolds would be products of a torus and a
certain $SO(m)$ Lie group.

We now introduce the controls in System (\ref{r588}) by the above
scheme. Using Equation (\ref{red_control}) and the fact that any
left-invariant vector field on the Lie group $SO(n)$ in our
redundant coordinates can be written as
$$
\dot{\cal X} = {\cal X} \Lambda, \quad \dot{\cal Y} = {\cal Y}
\Lambda, \quad \Lambda \in \mathfrak{so}(n),
$$
we get at once the following controlled system on the group:
\begin{eqnarray} \label{so_control}
\dot {\cal X}&={\cal X} [U ({\cal X} ^T {\cal Y} - {\cal Y} ^T {\cal
X})
+ {\cal X} ^T {\cal Y} U] - {\cal X} \left( \sum_i u_i \Lambda_i \right), \nonumber \\
\dot {\cal Y}&={\cal Y} [U ({\cal X} ^T {\cal Y} - {\cal Y} ^T {\cal
X}) - {\cal Y} ^T {\cal X} U] - {\cal Y} \left( \sum_i u_i \Lambda_i
\right).
\end{eqnarray}
This system describes an $n$-dimensional rigid body with ''symmetric
flywheels'', which is a direct generalization of the Liouville
problem of the rotation of a variable body \cite{Liou}.

\begin{proposition} \label{soc}
On can choose two vectors $\Lambda_1$ and $\Lambda_2$, such that for
any choice of the inertia operator $U$, one can transform any vortex
manifold to any other vortex manifold for any momentum in space
$M_s$, using the corresponding two control functions $u_1(t)$ and
$u_2(t)$.
\end{proposition}

{\it Proof.} First, we notice that System (\ref{r588}) preserves
volume in the phase space of the redundant variables $\cal{X},
\cal{Y}$ (this can be checked by the direct computation, but the
general result of the existence of an invariant measure for a
reduced system (\ref{v_lev}) or (\ref{v_prav}) with $\lambda = 0$
follows from \cite{Koz_b}). As the Lie algebra $\mathfrak{so}(n)$ is
semi-simple, controllability of System (\ref{so_control}) follows
from Corollary \ref{cor1}, Appendix \ref{Sec:AppB}. Proposition
\ref{soc} follows now from Theorem \ref{Th_vort}. $\Box$\\

\section{Conclusion and acknowledgements}

In this article, we considered the reduction of geodesic flows of
left- or right-invariant metrics on Lie groups to the group. The
reduced vector field has a remarkable hydrodynamic interpretation:
it is a velocity field for a stationary flow of an ideal fluid, the
the right- or left-invariant symmetry fields of the reduced field
being vortex vector fields, i.e., they annihilate the vorticity
2-form. The distribution of the vortex fields is always integrable,
thus it defines a manifold (at least locally), that we call a vortex
manifold. Typically, the vortex manifolds are tori.

We studied the following control problem. Consider a mechanical
system, whose configuration space is a Lie group and whose
Lagrangian is invariant to left translations on that group, and
assume that the mass geometry of the system may change under the
action of internal control forces. Such system can also be reduced
to the Lie group; with no controls, it describes a geodesic flow of
the left-invariant metric, given by the Lagrangian, and thus its
reduced flow is a stationary flow of an ideal fluid.

The control problem for such system is to find the conditions, under
which the system can be brought from any initial position in the
configuration space to another preassigned position by changing its
mass geometry. We showed that under these conditions, by changing
the mass geometry, one can also bring one vortex manifold to any
other preassigned vortex manifold. As an example, we considered the
$n$-dimensional Euler top. We wrote down the reduced controlled
system explicitly, showed that the vortex manifolds are tori, and
proved that, by changing the mass geometry, every such torus can be
transformed to any other torus.

The author wishes to thank Valery V. Kozlov for discussing the work.

\bigskip
\bigskip

\appendix{\Large \bf Appendix}

\section{Reduction to the Lie group for nonholonomic systems}

Consider the following equations, that we will refer to as the
generalized Euler equations (the left-invariant case):
 \be \label{EU1}
\dot {m_c} = ad^* _ {A(g) m_c} m_c,\\  \label{EU2} L_{g^{-1}*}
\dot{g} = A(g) m_c.
 \ee
Here $A(g): \mathfrak {g}^* \to \mathfrak {g}$ is positive definite
symmetric operator.

{\bf Example.} Consider the Chaplygin problem of a rigid ball
rolling on a horizontal plane. The equations of motion are:
 \bew
\dot{M} = M \times \omega, \quad \dot{\gamma} = \gamma \times
\omega, \\
M = I\omega + D\gamma \times (\omega \times \gamma),
 \eew
where $M$ is the ball momentum with respect to the moving axes,
fixed in the ball (momentum in the body), $\omega$ is the angular
velocity, $\gamma$ is the unit vertical vector, also written with
respect to the moving axes, matrix $I$ is the inertia tensor and $D$
is a constant. One can see that these equations are of the form
(\ref{EU1}).

\begin{theorem} \label{Thm_E} {\bf (The Euler theorem)}
The momentum in space $m_s$ is constant for the generalized Euler
equations (\ref{EU1}-\ref{EU2}).
\end{theorem}

{\it Proof.} Differentiate relation (\ref{mcs}) by time and apply
(\ref{EU2}), cf. \cite{AK}. $\Box$

Proposition \ref{Pr1} relied only on relation (\ref{mcs}), which
turns out to be true also for this case. Thus, reduction to the
group is possible, and the reduced vector field is
\begin{equation} \label{v_lev_A}
v (g) =L _ {g *} A(g) Ad^*_g m_s
\end{equation}
\begin{equation} \label{v_prav_A}
v (g) =R _ {g *} A(g) Ad^*_{g^{-1}} m_c.
\end{equation}
correspondingly in the left- or right-invariant case. Here $m_s$,
respectively $m_c$, is constant.

It would be interesting to find hydrodynamic description of the
reduced field.

\section{Controllability conditions} \label{Sec:AppB}

As system (\ref{red_control}) has the standard form
 \be \label{clas_contr}
\dot{x} = f(x) + \sum_i u_i g_i(x), \quad |u(t)| \le \epsilon
 \ee
of a classical control system, one can apply general theorems to it.

\begin{theorem} \label{Th_control_1}
Let the Lie group $G$ be compact. Then system (\ref{red_control}) is
controllable for all $\epsilon > 0$, if the minimal Lie subalgebra
of vector fields on $G$, which contains both vector fields $L_{g *}
\lambda_i$ and $v(g)$, spans the tangent space $T_g G$ at any $g \in
G$.
\end{theorem}

From Theorem \ref{Th_control_1} follows, that the minimal number of
controls is mainly defined by the Lie algebra structure. This
minimal number of controls should not necessarily be equal to the
number of the degrees of freedom of the system -- {\it unless the
Lie algebra is commutative}. If the dimension of the Lie algebra
$\mathfrak{g}$ is greater than 1, and we are interested in
controllability for {\it all} inertia operators and all momenta
$m_s$, then the minimal number of controls should necessarily be
greater than 1.

{\bf Example.} For controllability of a reduced Euler top on $SO(3)$
it is enough to have only one control, provided all the principal
axes of the inertia ellipsoid are different, and the momentum is the
space $m_s \ne 0$ is not directed along any of the principal axes --
see \cite{MK}. Obviously, if the ellipsoid of inertia is a sphere,
then the minimal number of controls is 2.

{\it Proof.} We only have to check that the vector field
$$
v(g) = L_{g *}  A^{-1} Ad^*_g m_s
$$
is Poisson-stable (i.e., almost all trajectories come back to the
vicinity of the initial conditions infinitely many times). Indeed,
if the minimal Lie subalgebra, which contains vectors $L_{g*}
\lambda_i$ and $v(g)$, spans the tangent space $T_g G$ at any $g \in
G$, then for each point $g_0 \in G$, the set of points $g(t, g_0,
u(t))$, accessible by the controls $u(t)$ for $0 < t < T$, form an
open set (the point $g_0$ itself may belong to the boundary). Under
the condition of the Poisson stability, these sets can be joined
together to get the necessary trajectory, see, e.g., \cite{Lian} for
details. The Poisson stability follows from the existence of a
smooth invariant measure of the reduced system $\dot{g} = v(g)$, as
we have assumed that the group $G$ is compact. But this is exactly
the case: if the Lie group $G$ is compact, the reduced system always
preserves a bi-invariant Haar measure on $G$, see \cite{Koz_b}.
$\Box$

\begin{corollary} \label{cor1}
Under conditions of Theorem \ref{Th_control_1}, let the Lie
algebra~$\mathfrak{g}$ be real and semi-simple. Then two controls is
sufficient for controllability for all values of the inertia
operator and of the momentum $m_s$.
\end{corollary}

{\bf Proof.} It is well known that a real semisimple Lie algebra is
generated by 2 elements, see, e.g., \cite{Bour}. $\Box$

{\bf Remark.} If the momentum $m_s = 0$, then under conditions of
Theorem \ref{Th_control_1}, system (\ref{red_control}) is
controllable even if the Lie group is noncompact -- this is the
classical Rashevsky-Chow theorem, see, e.g., \cite{Her}

The condition of Theorem \ref{Th_control_1} is usually referred to
as the Lie algebra rank condition (see, e.g., \cite{MK}). In real
systems, it may be difficult to check it directly, as, in principle,
the number of commutators one has to take is not bounded from above.
We suggest using a ''transversality'' condition, which is in the
next section.

\section{Transversality conditions}

The Lie algebra rank condition may be difficult to check, as the
number of commutators one should take to check it is not bounded
from below. We introduce a ''transversality'' condition, which seems
to be easier to verify in applications. This condition will also
provide stronger controllability results (in the non-analytic case):
we give a (rather trivial) example, when the Lie algebra rank
condition is not fulfilled, while the system is controllable.

Suppose that $N$ is a domain or a submanifold of $M$. We call a
vector field $v(x)$ {\it transversal} to $N$, if every phase
trajectory of the field $v(x)$, that starts in $N$ at $t = 0$,
escapes from $N$ both for some $t < 0$ and some $t > 0$. Notice that
it is enough to claim only one inequality (i.e., for example, say
that the trajectory leaves $N$ for some $t > 0$), if $v(x)$
preserves measure on $M$ and both $M$ and $N$ are compact: if at
least one inequality is fulfilled, there may not be any stationary
points of the field $v$ in $N$.

By a {\it finite system of commutators} we will understand a system
of the vector fields, which consists of the original fields $f(x)$,
$g_i(x)$ and a finite number of vector fields, obtained by taking
some fixed number of commutators of $f$, $g_i$, $[f,g_i]$, etc. If
the finite system of commutators has rank $n$ at least in one point
of $M$, then it has rank $n$ almost everywhere on $M$.

\begin{theorem} \label{Th_trans}
Suppose that all functions are analytic, and there exists a finite
system of commutators, which has rank $n$ on $M/N$, $N \in M$.
Suppose that the field $f(x)$ is transversal to $N$. Then the Lie
algebra rank condition is fulfilled.
\end{theorem}

{\it Proof.} Obviously the dimension of $N$ is less than $n$ -- due
to the analyticity condition. At any point $x_0 \in N$, take local
coordinates, such that $f(x) = (1, 0, \dots, 0)$, and $N$ is given
by condition $x_1 = 0$.

We give a proof in a simple situation, when the rank of a finite
system of commutators falls by 1 on $N$. At $x = x_0$, let $h_1,
\dots, h_{n-1}$ be independent vector fields, obtained as linear
combinations of $g_i$, $[f, g_i]$, etc., such that on $N$, $h_j^2 =
0$ for all $j$. Let the rank of the vector fields $f, h_2, \dots,
h_{n-1}$ be $n - 1$ at $x = x_0$. Then the rank of the vector fields
$f, [f, [f, \dots [f, h_1] \dots ]], h_2, \dots, h_{n-1}$ is $n$ at
$x = x_0$: otherwise, all the derivatives $\partial^k h_1^2/\partial
x_1^k = 0$, while $h_1^2 \ne 0$ for $x_1 \ne 0$. $\Box$

{\bf Example.} Consider the vector fields
 \bew
f(x) = (1, 0), g(x) = (0, 1 - \cos(x_1)); \quad x = (x_1
\pmod{2\pi}, x_2).
 \eew
These fields are independent only for $x_1 \ne 0$. The vector field
$f$ is transversal to the line $x_1 = 0$, and one can easily check
that the vector fields $f$ and $[f, [f,g]]$ are independent in the
neighbourhood of $x_1 = 0$.

In the non-analytic case, Theorem \ref{Th_trans} is not true (we
give a simple example below). However, a system may still be
controllable, even if the rank condition is not satisfied. One can
for example imagine a situation, when a system cannot be controlled
on some domain of the phase space.

\begin{theorem} \label{Th_c_t_1}
Let $M$ be compact, and let a finite system of commutators span the
tangent space $T_x M$ at any $x \in M/N$, $N \in M$. Suppose that
the vector field $f(x)$ is transversal to $N$. Then system
(\ref{clas_contr}) is controllable for all $\epsilon
> 0$.
\end{theorem}

{\it Proof.} Under the conditions of the theorem, we can reach every
point $x \in M/N$. For any point $x_0 \in N$, take a phase
trajectory that starts at $x_1 \in M/N$ and passes through $x_0$,
which exists due to the transversality condition. $\Box$

{\bf Example.} Consider the following system:
 \bew
\dot{x_1} = 1, \quad \dot{x_2} = u(t) g(x_1),
 \eew
where we take $x_1 \pmod{2\pi}$. The function $u(t)$ is a control,
which satisfies the condition $|u(t)| < \epsilon$, and
 \bew
g(x_1) = 0 \quad \mbox{for} \quad 0 < x_1 < \pi, \quad
            g(x_1) = \sin x_1 \quad \mbox{for} \quad  \pi < x_1 < 2\pi.
 \eew
The vector field $(1, 0)$ is transversal to the domain $0 \le x \le
\pi$, and is obviously Poisson-stable. This system is controllable
for all $\epsilon > 0$.

\end{document}